\documentclass[12pt,oneside]{amsart}
\usepackage[letterpaper,left=3cm,right=3cm,top=3cm,bottom=3cm]{geometry}

\usepackage{amssymb}
\newtheorem{theorem}{Theorem}[section]
\newtheorem{lemma}[theorem]{Lemma}

\newtheorem{e-definition}[theorem]{Definition\rm}
\newtheorem{rem}[theorem]{Remark\/}

\newcommand{\bN}{\mathbb{N}}
\newcommand{\bC}{\mathbb{C}}

\newcommand{\cB}{\mathcal{B}}
\newcommand{\cG}{\mathcal{G}}
\begin{document}
\title[An almost greedy uniformly bounded orthonormal basis]{An
  example of an almost greedy uniformly bounded orthonormal basis for
  $L_p([0,1])$}
\author{Morten Nielsen}
\address{Department of Mathematics\\ Washington University\\ 
Campus Box 1146\\ St.\ Louis, MO 63130, USA}
\email{mnielsen@math.wustl.edu}

\begin{abstract}
  We construct a uniformly bounded orthonormal almost greedy basis for
  $L_p([0,1])$, $1<p<\infty$. The example shows that it is not
  possible to extend Orlicz's theorem, stating that there are no
  uniformly bounded orthonormal unconditional bases for $L_p([0,1])$, $p\not=2$,
to the class of almost greedy bases.

\end{abstract}
\keywords{Bounded orthonormal systems, Schauder basis,
quasi-greedy basis, almost greedy basis, decreasing rearrangements}
\maketitle
\section{Introduction}
Let $\mathcal{B}=\{e_n\}_{n\in\bN}$ be a bounded Schauder
basis for a Banach space $X$, i.e., a basis for which
$0<\inf_n\|e_n\|_{X}\leq\sup_n\|e_n\|_X<\infty$. An approximation
algorithm associated with $\cB$ is a sequence $\{A_n\}_{n=1}^\infty$ of
(possibly nonlinear) maps $A_n:X\rightarrow X$ such that for $x\in X$,
$A_n(x)$ is a linear combination of at most $n$ elements from
$\cB$. We say that the algorithm is convergent if
$\lim_{n\rightarrow\infty} \|x-A_n(x)\|_X=0$ for every $x\in X$. For a
Schauder basis there is a natural convergent approximation
algorithm. Suppose the dual system to $\cB$ is given by
$\{e_k^*\}_{k\in\bN}$. Then the linear approximation algorithm is
given by the partial sums $S_n(x)=\sum_{k=1}^n e_k^*(x)e_k$.

Another quite natural approximation algorithm is the greedy
approximation algorithm where the partial sums are obtained by
thresholding the expansion coefficients.  The algorithm is defined as
follows.  For each element $x\in X$ we define the greedy ordering of
the coefficients as the map $\rho:\bN\rightarrow \bN$ with
$\rho(\bN)\supseteq \{j: e^*_j(x)\not=0\}$ such that for $j<k$ we have
either $|\,e^*_{\rho(k)}(x)|<|\,e^*_{\rho(j)}(x)|$ or
$|\,e^*_{\rho(k)}(x)|=|\,e^*_{\rho(j)}(x)|$ and $\rho(k)>\rho(j)$.
Then the greedy $m$-term approximant to $x$ is given by
$\mathcal{G}_m(x)=\sum_{j=1}^m e^*_{\rho(j)}(x)e_{\rho(j)}$.  The
question is whether the greedy algorithm is convergent.  This is
clearly the case for an {\em unconditional basis} where the expansion
$x=\sum_{k=1}^\infty e_k^*(x)e_k$ converges regardless of the
summation order. However, Temlyakov and Konyagin \cite{MR1716087}
showed that the greedy algorithm may also converge for certain
conditional bases. This leads to the definition of a quasi-greedy
basis.
\begin{e-definition}[\cite{MR1716087}]
  A bounded Schauder basis for a Banach space $X$ is called {\em
    quasi-greedy} if there exists a constant $C$ such that for $x\in
  X$, $\|\cG_m(x)\|_X\leq C\|x\|_X$ for $m\geq 1$.
\end{e-definition}
\noindent Wojtaszczyk proved the following result which gives a more 
 intuitive interpretation of quasi-greedy bases.

\begin{theorem}[\cite{MR1806955}]
  A bounded Schauder basis for a Banach space $X$ is {\em
    quasi-greedy} if and only if\, $\lim_{m\rightarrow \infty}\|x-\mathcal{G}_m(x)\|_X=0$ for
  every element $x\in X$.
\end{theorem}

In this note we study quasi-greedy bases for $L_p:=L_p([0,1])$,
$1<p<\infty$, with a particular structure. We are interested in
uniformly bounded bases $\mathcal{B}=\{e_n\}_{n\in\bN}$ such that
$\mathcal{B}$ is an orthonormal basis for $L_2$. It is a
well-known result by Orlicz that such a basis can be
unconditional only for $p=2$, so it is never trivially quasi-greedy except
for $p=2$. 

It was proved by Temlyakov \cite{MR1646563} that the trigonometric
system in $L_p$, $1\leq p\leq \infty$, $p\not=2$, fails to be
quasi-greedy.  Independently, and using a completely different
approach, C{\'o}rdoba and Fern{\'a}ndez \cite{MR1613739} proved the
same result in the range $1\leq p<2$.  One can also verify that the
Walsh system fails to be quasi greedy in $L_p$, $p\not=2$. This leads
to a natural question: are there any uniformly bounded orthonormal
quasi-greedy bases for $L_p$?

A negative answer to this question would give a nice improvement of
Orlicz's theorem to the class of quasi-greedy bases.  However, such an
improved result is not possible. Below we construct a uniformly
bounded orthonormal {\em almost greedy} basis for $L_p$,
$1<p<\infty$. An almost greedy basis is a quasi-greedy basis with one
additional property.
\begin{e-definition}
  A bounded Schauder basis $\{e_n\}_n$ for a Banach space $X$ is
  almost greedy if there is a constant $C$ such that for $x\in X$,
$$\|x - \mathcal{G}_n(x)\|_X\leq C \inf\big\{\big\|
x -\sum_{j\in A}
e^*_j(x)e_j\big\| : |A| = n, \alpha_j\in\bC, j\in A,n \in \bN\big\}.$$
\end{e-definition}
It was proved in \cite{MR1998906} that a basis is almost greedy if and
only if it is quasi-greedy and democratic. A Schauder basis
$\{e_n\}_n$ is called democratic if there exists $C$ such that for any finite
sets $A,B\subset\bN$ with $|A|=|B|$, we have
$$\big\|\sum_{j\in A} e_j\big\|_X\leq C\big\|\sum_{j\in B} e_j\big\|_X.$$

We can now state the main result of this note.

\begin{theorem}\label{thm}
  There exists a uniformly bounded orthonormal almost greedy basis for
  $L_p([0,1])$, $1<p<\infty$.
\end{theorem}
We should note that without the assumption that the system be
uniformly bounded, one can obtain a stronger result.  It is known that
the Haar system on $[0,1]$, normalized in $L_p$, is an unconditional
and democratic basis (a so-called greedy basis) for $L_p$,
$1<p<\infty$, see \cite{MR1628182}.

\section{A uniformly bounded almost greedy ONB for $L_p$}\label{s:3}
Classical uniformly bounded orthonormal system such as the
trigonometric basis and the Walsh system fail to form quasi-greedy
bases for $L_p:=L_p([0,1])$, $1<p<\infty$. The problem behind this
failure is that such systems are very far from being democratic.  This
behavior is not representative for all uniformly bounded orthonormal
bases. In this section we construct an example of a quasi greedy
uniformly bounded system in $L_p $, $1<p<\infty$. The example we
present is a variation on a construction by Kostyukovsky and
Olevski{\u\i} \cite{MR1618859}. The example in \cite{MR1618859} was
used to study pointwise convergence a.e.\ of greedy approximants to
$L_2$-functions. Later Wojtaszczyk \cite{MR1806955} used the same type
of construction to define quasi-greedy for $X\oplus\ell_2$, with $X$ a
quasi-Banach space with a Besselian basis.

Let us introduce some notation. The Rademacher functions are given by
$r_k(t)=\text{sign}(\sin(2^k\pi t))$ for $k\geq 1$. Khintchine's
inequality will be essential for the estimates below. The inequality
states that for $1\leq p<\infty$ there exist $A_p, B_p$ such that for
any finite sequence $\{a_k\}_{k\geq 1}$,
\begin{equation}
  \label{eq:khint}
 A_p \bigg(\sum_k |a_k|^2\bigg)^{1/2}\leq
\bigg(\int_0^1\bigg|\sum_k a_k r_k(t)\bigg|^p\,dt\bigg)^{1/p}\leq
B_p\bigg(\sum_k |a_k|^2\bigg)^{1/2}.
\end{equation}
Khintchine's inequality shows that the Rademacher functions form a
democratic system in $L_p$. However, the Rademacher system
is far from complete so it cannot be used directly to obtain an almost
greedy basis in $L_p$. In our example we use the fact
that the Rademacher functions form a subsystem of a complete system, namely the Walsh system.
The Walsh system $\mathcal{W}=\{W_n\}_{n=0}^\infty$ is defined as follows. For
$n=\sum_{j=1}^k \varepsilon_j 2^{j-1}$, the binary expansion of $n\in\bN$, we let
\begin{equation}\label{Walsh}
  W_n(t)=\prod_{j=1}^k r_j^{\varepsilon_j}(t).\end{equation}
The Walsh system forms a uniformly bounded orthonormal basis for $L_2$
and a 
Schauder basis for $L_p$, $1<p<\infty$, see \cite{MR1007141}. 
The idea is to reorder the Walsh system such that we obtain large dyadic blocks of Rademacher
functions with the remaining Walsh functions placed in between the
Rademacher blocks. Let us consider the details.

For $k=1,2,\ldots$, we define the $2^k\times 2^k$ Olevski{\u\i} matrix
$A^k=(a_{ij}^{(k)})_{i,j=1}^{2^k}$ by the following formulas
$$a_{i1}^{(k)}=2^{-k/2}\quad\text{for}\quad i=1,2,\ldots,2^k,$$
and for $j=2^s+\nu$, with $1\leq \nu\leq 2^s$ and $s=1,2,\ldots,k-1$, we let
$$a_{ij}^{(k)}=\begin{cases}2^{(s-k)/2}&\text{for}\quad (\nu-1)2^{k-s}<i\leq (2\nu-1)2^{k-s-1}\\
  -2^{(s-k)/2}&\text{for}\quad (2\nu-1)2^{k-s-1}<i\leq \nu 2^{k-s}\\
  0&\text{otherwise.}
\end{cases}$$ One can check (see \cite[Chapter IV]{MR0470599}) that
$A^k$ are orthogonal matrices and there exists a finite constant $C$
such that for all $i, k$ we have
\begin{equation}\label{e2}
\sum_{j=1}^{2^k} |a_{ij}^{(k)}|\leq C.
\end{equation}
Put $N_k=2^{10^k}$ and define $F_k$ such that $F_0=0$, $F_1=N_1-1$ and
$F_{k}-F_{k-1}=N_k-1$, $k=1,2,\ldots$.  We consider the Walsh
system $\mathcal{W}=\{W_n\}_{n=0}^\infty$ on $[0,1]$. We split
$\mathcal{W}$ into two subsystems. The first subsystem
$\mathcal{W}_1=\{r_k\}_{k=1}^\infty$ is the Rademacher functions with
their natural ordering. The second subsystem
$\mathcal{W}_2=\{\phi_k\}_{k=1}^\infty$ is the collection of Walsh
functions not in $\mathcal{W}_1$ with the ordering from
$\mathcal{W}$. We now impose the ordering
$$\phi_1,r_1,r_2,\ldots,r_{F_1},\phi_2,r_{F_1+1},\ldots, r_{F_2}, 
\phi_3,r_{F_2+1},\ldots,r_{F_3},\phi_4,\ldots$$
The block $\mathcal{B}_k:=\{\phi_k,r_{F_{k-1}+1},\ldots, r_{F_{k}}\}$
has length $N_k$, and we apply $A^{10^k}$ to $\mathcal{B}_k$ to obtain
a new orthonormal system $\{\psi_{i}^{(k)}\}_{i=1}^{{N_k}}$ given by
\begin{equation}\label{e1}
\psi_{i}^{(k)}=\frac{\phi_k}{\sqrt{N_k}}+\sum_{j=2}^{N_k}a_{ij}^{(10^k)}r_{F_{k-1}+j-1}.
\end{equation}
The system ordered $\psi^{(1)}_1,\ldots, \psi^{(1)}_{N_1},
\psi^{(2)}_1,\ldots, \psi^{(2)}_{N_2},\ldots$ will be denoted
$\mathcal{B}=\{\psi_k\}_{k=1}^\infty$. It is easy to verify that
$\mathcal{B}$ is an orthonormal basis for $L_2 $ since each matrix
$A^{10^k}$ is orthogonal. The system is uniformly bounded which
follows by \eqref{e2} and the fact that $\mathcal{W}$ is uniformly
bounded. The system $\mathcal{B}$ is our candidate for an almost
greedy basis for $L_p$, $1<p<\infty$.  We split the proof of Theorem
\ref{thm} into three parts. First we prove that $\mathcal{B}$ is
democratic in $L_p$. Then we prove that the system forms a Schauder
basis for $L_p$, and the final step is to prove that the system forms
a quasi-greedy basis for $L_p$.
\begin{lemma}\label{lem}
  The system $\mathcal{B}=\{\psi_k\}_{k=1}^\infty$ is democratic
  in $L_p $, $1<p<\infty$, with
$$\bigg\|\sum_{k\in A}\psi_k\bigg\|_p\asymp |A|^{1/2}.$$
\end{lemma}
\begin{proof}
  Fix $2<p<\infty$. Let $S=\sum_{k\in\Lambda} \psi_k$ with
  $|\Lambda|=N$. We write
$$S=\sum_{k=1}^\infty\sum_{j\in \Lambda_k} \psi^{(k)}_j=
\sum_{k=1}^\infty
\frac{|\Lambda_k|}{\sqrt{N_k}}\phi_k+\sum_{k=1}^\infty
\sum_{i\in\Lambda_k}
\sum_{j=2}^{N_k}a_{ij}^{(10^k)}r_{F_{k-1}+j-1}:=S_1+S_2,$$ with
$\sum_k|\Lambda_k|=|\Lambda|$, and $|\Lambda_k|\leq N_k$. Notice that
the coefficients of $\sum_{j\in \Lambda_k} \psi^{(k)}_j$ relative to
the block $\mathcal{B}_k$ has $l_2$-norm $|\Lambda_k|^{1/2}$ since
$A^{10^k}$ is orthogonal.  Hence, by Khintchine's inequality,
$$\|S_2\|_p=\bigg\|\sum_{k=1}^\infty
\sum_{j=2}^{N_k}\bigg(\sum_{i\in\Lambda_k}a_{ij}^{(10^k)}\bigg)r_{F_{k-1}+j-1}\bigg\|_p\leq
B_p \bigg(\sum_k |\Lambda_k|\bigg)^{1/2}=B_pN^{1/2}.$$ We now estimate
$S_1$. Write
$$S_1=\sum_{k=1}^\infty
\frac{|\Lambda_k|}{\sqrt{N_k}}\phi_k=\sum_{k\in A}
\frac{|\Lambda_k|}{\sqrt{N_k}}\phi_k+\sum_{k\in B}
\frac{|\Lambda_k|}{\sqrt{N_k}}\phi_k:=S_1^1+S_1^2,$$
where $A=\{k:|\Lambda_k|\leq (N_k)^{3/4}\}$ and
$B=\{k:|\Lambda_k|> (N_k)^{3/4}\}$. Using the Cauchy-Schwartz inequality,
$$\|S_1^1\|_p\leq
\sum_{k\in A} \frac{|\Lambda_k|}{\sqrt{N_k}}\leq
\bigg(\sum_{k\in A} |\Lambda_k|\bigg)^{1/2}
\bigg(\sum_{k\in A} \frac{|\Lambda_k|}{N_k}\bigg)^{1/2}\leq
N^{1/2}\bigg(\sum_{k\in \bN} N_k^{-1/4}\bigg)^{1/2}=CN^{1/2}.
$$
We turn to $S_1^2$. If $B$ is empty, we are done. Otherwise,
$B$ is a finite set and we can define $L=\max B$.
We have
$$
\sum_{k\in B;k<L}\frac{|\Lambda_k|}{\sqrt{N_k}}\leq 
\sum_{k\in B;k<L}\sqrt{N_k}\leq
\sum_{j=1}^{10^{L-1}}2^{j/2}\leq 2\cdot 2^{(10^{L-1}/2)}\leq
2\cdot 2^{(10^L/4)}\leq 2\frac{|\Lambda_L|}{\sqrt{N_L}}
$$
Hence,
$$
\|S_1^2\|_p\leq \sum_{k\in B}\frac{|\Lambda_k|}{\sqrt{N_k}}
\leq 3\frac{|\Lambda_L|}{\sqrt{N_L}}\leq
3\sqrt{|\Lambda_L|}\leq 3\sqrt{|\Lambda|}=3N^{1/2},
$$
where we used that $|\Lambda_L|\leq N_L$.  We conclude that
$\|S\|_p\leq C'N^{1/2}$, with $C'$ independent of $\Lambda$. Since
$N^{1/2}=\|S\|_2\leq \|S\|_p$ we deduce that $\mathcal{B}$ is
democratic in $L_p$, $2\leq p<\infty$. For $1<q<2$ we have
$\|S\|_q\leq \|S\|_2=N^{1/2}$. By H\"older's inequality,
$$N=\|S\|_2^2\leq \|S\|_q\|S\|_p\leq C_pN^{1/2}\|S\|_q,$$
for $1/q+1/p=1$. Again, we conclude that $\|S\|_q\asymp N^{1/2}$, so
 $\mathcal{B}$ is democratic in $L_q$, $1<q<2$.
\end{proof}
Next we prove that $\mathcal{B}$ is a basis for $L_p$. 
\begin{lemma}\label{lem2}
  The system $\mathcal{B}=\{\psi_k\}_{k=1}^\infty$ is a Schauder basis for $L_p $,
  $1<p<\infty$.
\end{lemma}

\begin{proof}
Notice that
  $\text{span}(\mathcal{B})=\text{span}(\mathcal{W})$ by
  construction, so $\text{span}(\mathcal{B})$ is dense in $L_p$,
  $1<p<\infty$, since $\mathcal{W}$ is a Schauder basis for $L_p$. 
Fix $2<p<\infty$ and let $f\in L_p$. Let
  $S_n(f)=\sum_{k=1}^n \langle f,\psi_k\rangle\psi_k$ be the partial
  sum operator. We need  to prove that the family of operators $\{S_n\}_n$ is
  uniformly bounded on $L_p$. Notice that $\{\langle f,\psi_k\rangle\}_k\in
  \ell_2(\bN)$ since $L_p\subset L_2$. For $n\in\bN$, we can find
  $L\geq 1$ and $1\leq m\leq N_L$ such that
$$S_n(f)=\sum_{k=1}^n \langle
f,\psi_k\rangle\psi_k=\sum_{k=1}^{L-1}\sum_{j=1}^{N_k} 
\langle f,\psi_j^{(k)}\rangle\psi_k^{(k)}+\sum_{k=1}^{m} 
\langle f,\psi_k^{(L)}\rangle\psi_k^{(L)}:=T_1+T_2.$$
Let us estimate $T_1$. If $L=1$ then $T_1=0$, so we may assume
$L>1$. The construction of $\mathcal{B}$ shows that $T_1$ is the
orthogonal projection of $f$ onto
$$\text{span}\bigg(\bigcup_{k=1}^{L-1}\bigcup_{j=1}^{N_k}\{\psi_k^{(k)}\}\bigg)=
\text{span}\big\{\{W_0,W_1,\ldots,
W_{L-2}\}\cup\{r_{\ell_0},r_{\ell_0+1},\ldots,
 r_{F_{L-1}}\}\big\},$$
with $\ell_0=\lfloor\log_2(L)\rfloor$.
It follows that we can rewrite $T_1$ as
$$T_1=\sum_{k=0}^{L-2}  \langle f,W_k\rangle W_k+P_R(f),$$
where $P_R(f)$ is the orthogonal projection of $f$ onto
$\text{span}\{r_{\ell_0},r_{\ell_0+1},\ldots, r_{F_{L-1}}\}$. Thus,
using Khintchine's inequality,
$$\|T_1\|_p\leq C_p\|f\|_p+B_p\|f\|_p,$$
where $C_p$ is the basis constant for the Walsh system in $L_p$. Next
we rewrite $T_2$ in the system $\{\phi_L,r_{F_{L-1}+1},\ldots,
r_{F_{L}}\}$,
$$T_2=\sum_{k=1}^{m}\langle f,\psi_k^{(L)}\rangle\frac{\phi_L}{\sqrt{N_L}}+ 
\sum_{j=2}^{N_L}\bigg(\sum_{k=1}^{m}
\langle f,\psi_k^{(L)}\rangle
a_{ij}^{(10^L)}\bigg)r_{F_{L-1}+j-1}.
$$
By Khintchine's inequality, and the fact that $A^{10^L}$ is orthogonal,
$$\bigg\|\sum_{j=2}^{N_L}\bigg(\sum_{k=1}^{m}
\langle f,\psi_k^{(L)}\rangle
a_{ij}^{(10^L)}\bigg)r_{F_{L-1}+j-1}\bigg\|_p\leq B_p\|\{\langle
f,\psi_k^{(L)}
\rangle\}_k\|_{\ell_2}<\infty.$$
Also,
$$\bigg\|\sum_{k=1}^{m}\langle
f,\psi_k^{(L)}\rangle\frac{\phi_L}{\sqrt{N_L}}\bigg\|_p \leq
\sum_{k=1}^{m}|\langle f,\psi_k^{(L)}\rangle|\frac{1}{\sqrt{N_L}} \leq
\|\{\langle f,\psi_k^{(L)}\rangle\}_k\|_{\ell_2}\sqrt{\frac m{N_L}}
\leq \|\{\langle f,\psi_k^{(L)}\rangle\}_k\|_{\ell_2},$$ so
$\|T_2\|_p<\infty$.  The estimates of $T_1$ and $T_2$ are independent
of $n$, and we obtain that $\sup_n \|S_n(f)\|_p<\infty$. Using the
Banach-Steinhaus theorem we deduce that $\{S_n\}_n$ is a uniformly
bounded family of linear operators on $L_p$. We conclude that
$\mathcal{B}$ is a Schauder basis for $L_p$, $2<p<\infty$, and the
result for $1<p<2$ follows by a duality argument.
\end{proof}

We can now complete the proof of Theorem \ref{thm}. Lemma \ref{lem3}
below together with Lemmas \ref{lem} and \ref{lem2} immediately give
Theorem \ref{thm}.
\begin{lemma}\label{lem3}
  The system $\mathcal{B}=\{\psi_k\}_{k=1}^\infty$ is a  quasi-greedy basis for $L_p $,
  $1<p<\infty$.
\end{lemma}
\begin{proof}
First we consider $2<p<\infty$. Let $f\in L_p\subset L_2$. Then we have the $L_p$-norm convergent expansion
\begin{equation}\label{eq:2}
f=\sum_{i=1}^\infty\langle f,\psi_i\rangle\psi_i,
\end{equation}
with $\|\{\langle f,\psi_i\rangle\}_i\|_{\ell_2}\leq
\|f\|_2\leq\|f\|_p$.  It suffices to prove that $\mathcal{G}_m(f)$ is
convergent in $L_p$ since $\mathcal{G}_m(f)\rightarrow f$ in $L_2$.
We write (formally)
\begin{align*}
  f&=\sum_{k=1}^{\infty}\sum_{j=1}^{N_k} \langle f,\psi_j^{(k)}\rangle\psi_k^{(k)}\\
  &=\sum_{k=1}^{\infty}\sum_{j=1}^{N_k}\langle
  f,\psi_j^{(k)}\rangle\frac{
    \phi_k}{\sqrt{N_k}}+\sum_{k=1}^\infty\sum_{i=1}^{N_k} \langle
  f,\psi_i^{(k)}\rangle
  \sum_{j=2}^{N_k}a_{ij}^{(10^k)}r_{F_{k-1}+j-1}\\
  &=S^1+S^2.
\end{align*}
Consider a sequence $\{\varepsilon^k_i\}\subset\{0,1\}$. By
Khintchine's inequality, and the fact that each $A^{10^k}$ is
orthogonal,
$$\bigg\|
\sum_{k=1}^\infty\sum_{j=2}^{N_k} \bigg(\sum_{i=1}^{N_k}
\varepsilon_i^k\langle f,\psi_i^{(k)}\rangle
a_{ij}^{(10^k)}\bigg)r_{F_{k-1}+j-1} \bigg\|_p\leq B_p
\bigg(\sum_k\sum_{i=1}^{N_k} \varepsilon_i^k|\langle
f,\psi_i^{(k)}\rangle|^2\bigg)^{1/2}.$$   It
follows that $S^2$ is convergent and actually converges
unconditionally in $L_p$. From this and the convergence of the series
\eqref{eq:2}, we conclude that the partial sums for the series $S^1$,
$$S^1_n=\sum_{k=1}^{L-1}\sum_{j=1}^{N_k}\langle f,\psi_j^{(k)}\rangle\frac{ \phi_k}{\sqrt{N_k}}+
\sum_{j=1}^{m}\langle f,\psi_j^{(L)}\rangle\frac{ \phi_L}{\sqrt{N_L}}
$$
converge in $L_p$. 

The series defining $S^2$ converges unconditionally, so it suffices to
prove that the series defining $S^1$ converges in $L_p$ when the
coefficients $\{\langle f,\psi_i\rangle\}_i$ are arranged in
decreasing order.  We define the sets
 \begin{align}
   \Lambda_k&=\bigg\{j:\frac1{N_k}< |\langle f,\psi_j^{(k)}
\rangle|<\frac1{N_k^{1/10}}\bigg\}\nonumber\\
   \Lambda_k'&=\bigg\{j: |\langle f,\psi_j^{(k)}
\rangle|\leq \frac1{N_k}\bigg\}\label{eq:lamb}\\
   \Lambda_k''&=\bigg\{j: |\langle f,\psi_j^{(k)}\rangle|\geq
   \frac1{N_k^{1/10}}\bigg\}\nonumber.
 \end{align}
Then (formally)
$$S^1=
\sum_{k=1}^{\infty}\sum_{j\in \Lambda_k}\langle
f,\psi_j^{(k)}\rangle\frac{ \phi_k}{\sqrt{N_k}} +
\sum_{k=1}^{\infty}\sum_{j\in \Lambda_k'}\langle
f,\psi_j^{(k)}\rangle\frac{ \phi_k}{\sqrt{N_k}} +
\sum_{k=1}^{\infty}\sum_{j\in \Lambda_k''}\langle
f,\psi_j^{(k)}\rangle\frac{ \phi_k}{\sqrt{N_k}}=T+T'+T''.$$ Notice
that $\sum_{j\in \Lambda_k'}\|\langle f,\psi_j^{(k)}\rangle\frac{
  \phi_k}{\sqrt{N_k}}\|_p\leq \sum_{j\in \Lambda_k'}\frac{|\langle
  f,\psi_j^{(k)}\rangle|}{\sqrt{N_k}}\leq 1/\sqrt{N_k}$, so the series
defining $T'$ converges absolutely in $L_p$. For $T''$ we notice that
$|\Lambda_k''|\leq \|f\|_2^2N_k^{1/5}$, so
$$\sum_{j\in \Lambda_k''}\bigg\|\langle f,\psi_j^{(k)}\rangle\frac{
  \phi_k}{\sqrt{N_k}}
\bigg\|_p\leq \sum_{j\in \Lambda_k''}\frac{|\langle f,\psi_j^{(k)}\rangle|}{\sqrt{N_k}}\leq 
\frac{\|f\|_p\|f\|_2^2}{{N_k^{3/10}}},$$ and the series defining $T''$
converges absolutely in $L_p$.

The series defining $S^1$, $T'$ and $T''$ converge in $L_p$, so we may
conclude that the series defining $T$ converges in $L_p$. From
\eqref{eq:lamb}, we get
$$|\langle f,\psi_i^{(k)}\rangle|>\frac1{N_k}\geq \frac1{N_{k+1}^{1/10}}\geq
|\langle f,\psi_j^{(k+1)}\rangle|,\qquad i\in\Lambda_k,j\in\Lambda_{k+1}; k=1,2,\ldots$$
so when we arrange $T$ by decreasing order, the rearrangement can only
take place 
inside the blocks. The estimate
$$\sum_{j\in \Lambda_k}\bigg\|\langle f,\psi_j^{(k)}\rangle\frac{
  \phi_k}{\sqrt{N_k}}\bigg\|_p
\leq \bigg(\sum_{j\in \Lambda_k} |\langle
f,\psi_j^{(k)}\rangle|^2\bigg)^{1/2}
\frac{|\Lambda_k|^{1/2}}{\sqrt{N_k}},\quad k\geq 1,$$
shows that rearrangements inside blocks are well-behaved, and
$$\sum_{j\in \Lambda_k}\bigg\|\langle f,\psi_j^{(k)}\rangle\frac{
  \phi_k}{\sqrt{N_k}}
\bigg\|_p\rightarrow 0\quad \text{as}\quad k\rightarrow \infty.$$
We conclude that $\mathcal{G}_m(f)$ is convergent in $L_p$ and
consequently $\cB$ is a quasi-greedy basis in $L_p$, $2\leq
p<\infty$. Fix $1<q<2$ and let $p$ be given by $1/q+1/p=1$. By Lemma
\ref{lem}, for any finite subset $A\subset\bN$,
$$\bigg\|\sum_{k\in A}\psi_k\bigg\|_q
\bigg\|\sum_{k\in A}\psi_k\bigg\|_p\leq C|A|,$$ so $\cB$ is a
so-called bi-democratic system in $L_p$. It follows from \cite[Theorem
5.4; $(1)\Rightarrow (2)$]{MR1998906} that $\cB$ is a quasi-greedy
basis for $L_q$. This completes the proof.
\end{proof}

\begin{rem}
  To get a uniformly bounded quasi-greedy basis consisting of smooth
  functions, we can use the same construction based on the
  trigonometric system with any lacunary subsequence playing the role
  of the Rademacher system.
\end{rem}

\subsection*{Acknowledgment}
The author would like to thank Boris Kashin and R\'emi Gribonval for
discussions on uniformly bounded quasi-greedy systems that provided
inspiration for our example.

\end{document}